\documentclass[11pt]{article}
\usepackage{graphicx}
\usepackage{amsmath}
\usepackage{amssymb}
\usepackage[utf8]{inputenc}
\usepackage[mathscr]{eucal}
\usepackage{fancyhdr}
\usepackage{fancybox}
\usepackage{color}
\usepackage{longtable}
\usepackage{dsfont}

\newtheorem {Proposition}{Proposition} [section]

\newtheorem {Theorem}[Proposition]{Theorem}

\newtheorem {Remark}{Remark}[Proposition]

\title{The statistical effect of entropic regularization in optimal transportation}
\author{Eustasio del Barrio and Jean-Michel Loubes}

\begin{document}

\maketitle
\begin{abstract}
We propose to tackle the problem of understanding the effect of regularization in Sinkhorn algotihms. In the case of Gaussian distributions we provide a closed form for the regularized optimal transport which enables to provide a better understanding of the effect of the regularization from a statistical framework.
\end{abstract}

\section{Introduction.} 

Statistical methods based on optimal transportation have received a considerable amount of attention in recent times. While the topic has a long history,
computational limitations (and also lack of a well developed distributional theory) hampered its applicability for years. Some recent advences (see Cuturi, Peyre, Schmitzer, Rigollet,...) have completely changed the scene and now statistical methods based on optimal transportation are everywhere (see, e.g.\cite{peyre2019computational},for Kernel based methods \cite{Kolouri_2016_CVPR,bachoc2017gaussian}, in Fair Machine Learning \cite{gordaliza2019obtaining}).

Monge-Kantorovich distances are defined using a cost function $c$ as
$$W_c(P,Q)= \min_{\pi \in \Pi(P,Q)} \int c(x,y) d\pi(x,y),$$
where $\Pi(P,Q)$ denotes the set of distributions  with marginals $P$ and $Q$. \\ Computing such distances requires to solve  in the discrete case a linear program. Actually solving the original discrete optimal transport problem, for two discrete distributions
$P=\sum_{i=1}^n p_i \delta_{x_i}$ and $Q=\sum_{i=1}^m q_i \delta_{y_i}$ and a cost matrix $c$ $c_{ij}=c(x_i,y_j)$ for all $(i,j)\in [1,n]\times [1,m]$,  amounts to solve the minimization with respect to $\pi$ the transportation plan \begin{equation} \label{otemp}\min_{ \pi \in \Pi(P,Q)} <c,\pi>  \end{equation}
 $\Pi(P,Q)=\{ \pi \in \mathbb{R}^{n\times n }_{+}, \: \pi \mathds{1}_n = P , \:  \pi^T \mathds{1}_m = Q \}$
where $ \pi \mathds{1}_n=(\sum_{j=1}^n \pi_{ij})_i$ and  $ \pi^T \mathds{1}_m=(\sum_{j=1}^m \pi_{ij})_j$. This minimization  is yet a linear problem (see \cite{kantorovich1942translocation}) but it turns to be  computationally difficult. Different algorithms have been proposed such as the Hungarian algorithm \cite{kuhn1955hungarian}, the simplex algorithm \cite{luenberger1984linear} or others versions using interior points algorithms \cite{orlin1988faster}. The complexity of these methods is at most of order $O(n^3 \log(n))$ for the OT problem between two discrete distributions with equal size $n$. \vskip .1in
To overcome this issue, regularization methods have been proposed to approximate the optimal transport problem by adding a penalty. The seminal paper by \cite{cuturi2013sinkhorn} provides the description of the Sinkhorn algorithm to regularize  optimal transport by using the  entropy of the transportation plan 
$ H(\pi)= \sum_{i,j} \pi_{ij} \log(\pi_{ij}),$ and changing the initial optimization program \eqref{otemp} into a strictly convex one
\begin{equation} \label{otempreg}\min_{ \pi \in \Pi(P,Q)} \{  <c,\pi>  + \varepsilon H(\pi) \}. \end{equation} 
The minimization of this criterion is achieved using Sinkhorn algorithm. We refer to \cite{peyre2019computational} and references therein for more details.
The introdution of the Sinkhorn divergence enables to obtain an $\varepsilon$-approximation of the optimal transport distance which can be computed, as pointed out in \cite{altschuler2017near},  with a  complexity of algorithm  of order $O(\frac{n^2}{\varepsilon^3})$, hence in a much faster way than the original optimal transport problem. Several toolboxes have been developed to compute regularized OT such among others as \cite{flamary2017pot} for python, \cite{klatt2017package} for R. \\
Other algorithms can be used to minimize \eqref{otempreg}. In~\cite{genevay2016stochastic} stochastic gradient
descent is applied to solve the entropy-regularized OT problem while in~\cite{dvurechensky2018computational} an accelerated gradient descent is proposed improving the complexity to  $O(\min (\frac{n^2}{\varepsilon^2}, \frac{n^{9/4}}{\varepsilon}))$.  \vskip .1in
The influence of the penalty is balanced introducing a parameter $\varepsilon >0 $ which controls the balance between the approximation of the optimal transport distance and its computational feasibility. Note also that others regularizing penalty have been proposed, for instance the entropy with respect to the product of marginals . \vskip .1in
 Beyond computational convenience, regularization has a statistical impact  and few results exist in the literature.  An enjoyable property of regularized optimal transport is that the convergence of its empirical version is faster than the standard optimal transport. Actually, if $P_n$ and $Q_n$ are empirical versions of distributions $P$ and $Q$ in $\mathbb{R}^d$, then Monge-Kantorovich distances suffer from the curse of dimensionality and converge under some assumptions at a rate at most $n^{-1/d}$, this rate may be improved under some assumptions as pointed out in \cite{weed2019sharp} for instance. As shown first in \cite{genevay2019sample} for  distributions defined on a bounded domain,  and sharpened for sub-Gaussian distributions in \cite{mena2019statistical}, the rate of convergence of regularized OT divergences is of order $\frac{1}{\varepsilon^{2+[5d/4]}} / {\sqrt{n}} $. \vskip .1in
In the recent years, optimal transport theory has been extensively used in unsupervised learning in order to characterize the mean of observations, giving rise to the notion of Wasserstein barycenters. This point of view is closely related to the notion of Fr\'echet means which has been used in statistics in preliminar works such as \cite{dupuy2011non} . The problem of existence and uniqueness of the Wasserstein barycenter of distributions $P_1,\dots,P_k$, where at least one of these distributions has a density, has been tackled in \cite{agueh2011barycenters}. The asymptotic property of Wasserstein barycenters have been studied in \cite{boissard2015distribution} or \cite{le2017existence}.
     However their computation is a difficult issue apart from the scatter-location family case. In this case a fixed point solution method can be derived to compute their barycenter as explained in \cite{alvarez2016fixed}. Hence, some authors have replaced Monge-Kantorovich distance by the Sinkhorn divergence and thus have considered the notion of Sinkhorn barycenter as  in \cite{cuturi2014fast} or \cite{bigot2019penalization}. In this setting, the distributions are discretized and the usual Sinkhorn's algorithm for discrete distributions is applied. Results proving the consistency of empirical Sinkhorn barycenters towards population Sinkhorn barycenters can be derived and the rate of convergence can be upper bounded by a bound depending on the number of observations, the discretization scheme and the trade-off parameter $\varepsilon$. Here again little is said to derive the statistical properties of the Sinkhorn barycenter and its property with respect to the original Wasserstein barycenter. \vskip .1in
 Hence, for both computational and statistical properties, the influence of $\varepsilon$ is crucial  and the results dealing with the approximation properties of regularized OT with regards to standard OT are scarce.  Very recently some papers all independently have proved in~\cite{janati2020entropic} similar expressions for the closed form of regularized optimal transport between Gaussian distributions, including in their case the case of unbalanced transport. In \cite{mallasto2020entropyregularized},  similar formulations have been derived using proofs based on the solution of the Schr\"odinger system that can be written to compute the entropic transport plan. All three points of view are complementary and provide new insights on entropic optimal transport.
 \vskip .1in
 Our contribution is the following

%Plan: complete the Intro expanding these items:
\begin{itemize}
% \item Describe entropic regularization and cite Cuturi's Sinkhorn algorithm
 % \item Mention the computation advantages coming from the regularization
 %\item Comment on regularized barycenters (cite Bigot et al papers)
 \item We investigate in this paper this impact
 \item We find that optimal regularized coupling of Gaussian measures is Gaussian and compute regularized transportation cost between Gaussians (Theorem \ref{GaussianEntropicTC})
 \item The Gaussian case is not just an interesting benchmark. In fact, just as in the classical (unregularized) optimal transportation problem,
 for probabilities with given means and covariance matrices the entropic transportation cost is minimized for Gaussian distributions. This is a generalization of Gelbrich lower bound the entropic setup (Theorem \ref{EntropicGelbrich}).
 \item Also as in the classical case, the entropic barycenter of Gaussian probabilities is Gaussian (Theorem \ref{entropic_Gaussian_barycenter}).
 \item Entropic variation around barycenter lower bounded by explicit expression from Gaussian case
 \item We see that entropic regularization basically amounts to smoothing via convolution with a Gaussian kernel, which results in
 added variance. The regularization parameter controls the increase in variance
 
% \item TODO: can this yield good guidelines to choose a compromise between computational time (better for large $\varepsilon$) and
% statistical accuracy (better for small $\varepsilon$)?
\end{itemize}

\section{Regularized optimal transport.}
 
We consider the entropic regularization of the transportation cost, namely, for 
probabilities $P$, $Q$ on $\mathbb{R}^d$, 
$$\mathcal{W}^2_{2,\varepsilon}(P,Q)=\min_{\pi\in\Pi(P,Q)} I_\varepsilon[\pi]$$ 
with
\begin{equation}\label{W2e}
I_\varepsilon[\pi]=\int_{\mathbb{R}^d\times \mathbb{R}^d} \|x-y\|^2 d\pi(x,y) 
+\varepsilon H(\pi).
\end{equation}
Here $H$ stands for the negative of the differential or  Boltzmann-Shannon entropy, that is, if
$\pi$ has density $r$ with respect to Lebesque measure on $\mathbb{R}^d\times \mathbb{R}^d$, then
$$H(\pi)=\int_{\mathbb{R}^d\times \mathbb{R}^d} r(x,y)\log r(x,y)dxdy,$$
while $H(\pi)=+\infty$ if $\pi$ does not have a density.

\medskip
The entropy term $H$ modifies the linear term in classical optimal transportation (the quadratic transportation cost)
to produce a strictly convex functional. This is not the only possible choice. Alternatively, we could
fix two reference probability measures on $\mathbb{R}^d$, say $\mu$ and $\nu$, and 
consider 
$$\mathcal{W}^2_{2,\varepsilon,\mu,\nu}(P,Q)=\min_{\pi\in\Pi(P,Q)} I_{\varepsilon,\mu,\nu}[\pi]$$ 
where
\begin{equation}\label{W2emunu}
I_{\varepsilon,\mu,\nu}[\pi]=\int_{\mathbb{R}^d\times \mathbb{R}^d} \|x-y\|^2 d\pi(x,y) 
+\varepsilon K(\pi|\mu\otimes \nu)
\end{equation}
and $K$ denotes the Kullback-Leibler divergence, namely, for probability measures, $\rho$, $\eta$, $K(\rho\|\eta)=\int \log \frac{d\rho}{d\eta} d\rho$
if $\rho\ll \eta$ and $K(\rho\|\eta)=+\infty$ otherwise. In the case when $\mu=\nu$ 
is the centered normal distribution on $\mathbb{R}^d$ with covariance matrix, $\lambda I_d$, for some $\lambda>0$, we will simply write
$I_{\varepsilon,\lambda}[\pi]$ and $\mathcal{W}^2_{2,\varepsilon,\lambda}(P,Q)$.

In our definitions of the regularized transportation cost we have written $\min$ instead of $\inf$. The existence of the minimizer follows easily.
In the case of $\mathcal{W}^2_{2,\varepsilon,\mu,\nu}(P,Q)$, for instance, let us assume that $\pi_n\in\Pi(P,Q)$ is a minimizing sequence, that is,
$I_{\varepsilon,\mu,\nu}[\pi_n]\to\inf_{\pi\in\Pi(P,Q)} I_{\varepsilon,\mu,\nu}[\pi]=m<\infty$. Since $\pi_n$ have fixed marginals, $\{\pi_n\}$
is a tight sequence and we can extract a weakly convergent subsequence, that we keep denoting $\pi_n$, say $\pi_n\to\pi_0$. Obviously $\pi_0\in \Pi(P,Q)$.
By Fatou's Lemma $\int_{\mathbb{R}^d\times \mathbb{R}^d} \|x-y\|^2 d\pi(x,y)\leq \liminf_n \int_{\mathbb{R}^d\times \mathbb{R}^d} \|x-y\|^2 d\pi_n(x,y)$
and by lower semicontinuity of relative entropy (see, e.g., Lemma 1.4.3 in \cite{DupuisEllis}) $K(\pi_0|\mu\otimes\nu)\leq \liminf_n K(\pi_n|\mu\otimes\nu)$.
But this shows that $I_{\varepsilon,\mu,\nu}[\pi_0]\leq \liminf_n I_{\varepsilon,\mu,\nu}[\pi_n]=m$, hence, $\pi_0$ is a minimizer. The case of $\mathcal{W}^2_{2,\varepsilon}(P,Q)$ follows similarly. Futhermore, if the transportation cost is finite then the minimizer is unique, since the relative entropy is strictly convex in its domain.

The choice of the reference measures is arbitrary. However, its influence on the regularized
optimal transport is limited. In fact, if we replace $\mu$, $\nu$ with equivalent measures $\mu'$, $\nu'$ (in the sense of $\mu$ and $\mu'$ being
mutually absolutely continuous with respect to each other and similarly for $\nu$ and $\nu'$) then $\pi\ll \mu\otimes \nu$ if and only if $\pi\ll \mu'\otimes \nu'$ and 
then $\frac{d\pi}{d(\mu'\otimes \nu')}=\frac{d\pi}{d(\mu\otimes \nu)}\big/ \Big(\frac {d\mu'}{d\mu} \frac {d\nu'}{d\nu} \Big)$.
Hence, for any $\pi\in \Pi(P,Q)$ with $\pi\ll \mu\otimes \nu$, writing $r=\frac{d\pi}{d(\mu\otimes \nu)}$ we have
\begin{equation}\label{red1}
K(\pi\| \mu\otimes \nu)-K(\pi\| \mu'\otimes \nu')=\int_{\mathbb{R}^d} \log\big(\textstyle \frac {d\mu'}{d\mu}(x)\big)dP(x)+\int_{\mathbb{R}^d} \log\big(\textstyle \frac {d\nu'}{d\nu}(y)\big)dQ(y)
\end{equation}
and we see that the difference does not depend on $\pi$. In particular the minimizer, if it exists, does not depend on the choice of $\mu$, $\nu$.
Furthermore, if $\mu$ and $\nu$ have a positive density on $\mathbb{R}^d$ then $I_\varepsilon[\pi]$ and $I_{\varepsilon,\mu,\nu}[\pi]$
differ only in a constant and, again, the minimizer does not depend on the choice of $\mu$, $\nu$. The minimal value, however, does depend
on the choice of the regularization term and this has an impact, for instance, in the barycenter problem, as we will see later.

\bigskip
We prove in this section that the entropic regularization of the transportation problem between nondegenerate Gaussian laws admits a
(unique) minimizer which is also Gaussian (on the product space). We provide a explicit expression for the mean and covariance
of this minimizer. Our proof is self-contained in the sense that we prove the existence of a minimizer in this setup. This existence 
could be obtained from more general results (see, e.g., Theorem 3.2 in \cite{Chizatetal} or Remark 4.19 in \cite{peyre2019computational}) based
on duality. We obtain the minimizer, instead, from the analysis of a particular type of matrix equation: the so-called algebraic Riccatti
equation. This equation has been extensively studied (see \cite{LancasterRodman1995}) and efficient numerical methods for the computation
of solutions are available (see, e.g., \cite{BiniIanazzoMeini2012}). However, the particular Riccatti equation which is of interest for
the entropic transportation problem (see (\ref{Riccatti})) has a particularly simple structure and
its unique positive definite solution admits an explicit expression. This is shown in our next result.

\begin{Proposition} \label{RiccattiSolutions} If $\Sigma_1$, $\Sigma_2$ are real, symmetric, positive definite $d\times d$ matrices and $\varepsilon > 0$ then the
unique symmetric, positive definite solution of the matrix equation
\begin{equation}\label{Riccatti}
X\Sigma_1X +\textstyle \frac \varepsilon  2 X= \Sigma_2 
\end{equation}
is
\begin{equation}\label{RiccattiSol}
X_\varepsilon=\Sigma_1^{-1/2} \big(\Sigma_1^{1/2}\Sigma_2 \Sigma_1^{1/2}+\textstyle (\frac \varepsilon 4)^2 I_d \big)^{1/2} \Sigma_1^{-1/2}-\textstyle \frac \varepsilon 4 \Sigma_1^{-1}.
\end{equation}
Furthermore, if 
$$\Sigma_\varepsilon=\left[
\begin{matrix}
\Sigma_1 & \Sigma_1 X_\varepsilon \\
X_\varepsilon \Sigma_1 &\Sigma_2
\end{matrix}
\right]
$$
then $\Sigma_\varepsilon$ is a real, symmetric, positive definite $2d\times 2d$ matrix and
$$\Sigma_\varepsilon^{-1}=\left[
\begin{matrix}
\Sigma_1^{-1}+\frac 2 \varepsilon X_\varepsilon & -\frac{2}{\varepsilon} I_d \\
-\frac{2}{\varepsilon} I_d &\frac 2 \varepsilon X_\varepsilon^{-1}
\end{matrix}
\right].
$$
\end{Proposition}

\bigskip
\noindent
\textbf{Proof.} The fact that $X_\varepsilon$ solves (\ref{Riccatti}) can be checked by simple inspection. $X_\varepsilon$ is obviously symmetric. Hence, it suffices to show that it is positive definite or, equivalently, that $\big(\Sigma_1^{1/2}\Sigma_2 \Sigma_1^{1/2}+\textstyle (\frac \varepsilon 4)^2 I_d \big)^{1/2}-\textstyle \frac \varepsilon 4 I_d$ is positive definite. 
This, in turn, will follow if we prove that every eigenvalue, 
say $\lambda$, of $\big(\Sigma_1^{1/2}\Sigma_2 \Sigma_1^{1/2}+\textstyle (\frac \varepsilon 4)^2 I_d \big)^{1/2}$ satisfies $\lambda >\frac\varepsilon 4$. 
But this is a consequence of the fact that the eigenvalues of
$\big(\Sigma_1^{1/2}\Sigma_2 \Sigma_1^{1/2}+\textstyle (\frac \varepsilon 4)^2 I_d \big)^{1/2}$ are $\sqrt{s+(\frac \varepsilon 4)^2}$ with $s$ ranging in the set of eigenvalues
of $\Sigma_1^{1/2}\Sigma_2 \Sigma_1^{1/2}$, which is positive definite. Consequently, $X_\varepsilon$ is a positive definite solution of (\ref{Riccatti}).
To prove uniqueness we set $Z=\frac \varepsilon 4I_d+\Sigma_1 X$ and note that if $X$ is a solution to (\ref{Riccatti}) then 
\begin{equation}\label{uniqueness}
XZ=\Sigma_2-\textstyle\frac \varepsilon 4 X. 
\end{equation}
But then $X=\Sigma_1^{-1}(Z-\frac \varepsilon 4 I_d)$ and substitution in (\ref{uniqueness}) yields
$$\Sigma_2+\textstyle\big(\frac\varepsilon 4 \big)^2\Sigma_1^{-1} =\Sigma_1^{-1} Z^2$$
or, equivalently,
$$Z^2=\Sigma_1\Sigma_2 +\textstyle\big(\frac\varepsilon 4 \big)^2 I_d.$$
Observe now that $A:=\Sigma_1^{-1/2}Z\Sigma_1^{1/2}$ is a symmetric, positive definite matrix. From the last identity
we see that
$$A^2=\Sigma_1^{-1/2}Z^2\Sigma_1^{1/2}=\Sigma_1^{1/2}\Sigma_2\Sigma_1^{1/2} +\textstyle\big(\frac\varepsilon 4 \big)^2 I_d.$$
Therefore, $A=(\Sigma_1^{1/2}\Sigma_2\Sigma_1^{1/2} +\textstyle (\frac\varepsilon 4 )^2 I_d)^{1/2}$. We conclude that, necessarily, $X=X_\varepsilon$.

We show next that $\Sigma_\varepsilon$ is positive definite. In fact, (see, e.g., Theorem 1.3.3 in \cite{Bhatia}) it suffices to show that
$\Sigma_1-\Sigma_1 X_\varepsilon \Sigma_2^{-1}X_\varepsilon \Sigma_1$ is positive definite. Since $X_\varepsilon$ solves (\ref{Riccatti}),
we have that $X_\varepsilon^{-1}\Sigma_2 X_\varepsilon^{-1}=\Sigma_1+\frac \varepsilon 2  X_\varepsilon^{-1}$ and the last condition becomes
that $\Sigma_1-\Sigma_1  (\Sigma_1+\frac  \varepsilon 2 X_\varepsilon^{-1})^{-1}\Sigma_1$ has to be positive definite. But this holds
if and only if
$$U=\left[
\begin{matrix}
\Sigma_1 & \Sigma_1 \\
 \Sigma_1 &\Sigma_1+ \frac \varepsilon 2 X_\varepsilon^{-1}
\end{matrix}
\right]
$$
is positive definite. Since $[x^T y^T]U[x^T y^T]^T=(x+y)^T\Sigma_1(x+y)+ \frac \varepsilon 2 y^T X_\varepsilon^{-1}y$, we conclude
that $\Sigma_\varepsilon$ is indeed positive definite.

To complete the proof we note that from the well known
identity for the inverse of block partitioned matrices
$$\Sigma_\varepsilon^{-1}=\left[ 
\begin{matrix}
(\Sigma_1-\Sigma_1 X_\varepsilon \Sigma_2^{-1} X_\varepsilon\Sigma_1 )^{-1} & -X_\varepsilon(\Sigma_2-X_\varepsilon \Sigma_1 X_\varepsilon )^{-1} \\
-(\Sigma_2-X_\varepsilon \Sigma_1 X_\varepsilon )^{-1}X_\varepsilon &(\Sigma_2-X_\varepsilon \Sigma_1 X_\varepsilon )^{-1} 
\end{matrix}
\right],
$$
Since $X_\varepsilon$ solves (\ref{Riccatti}) we have that $(\Sigma_2-X_\varepsilon \Sigma_1 X_\varepsilon )^{-1}=\frac 2 \varepsilon X_\varepsilon^{-1}$.
We similarly check that $(\Sigma_1-\Sigma_1 X_\varepsilon \Sigma_2^{-1} X_\varepsilon\Sigma_1 )(\Sigma_1^{-1}+\frac 2 \varepsilon X_\varepsilon)
={I_d+\frac 2 \varepsilon \Sigma_1 X_\varepsilon-\Sigma_1 X_\varepsilon \Sigma_2^{-1} X_\varepsilon-\frac 2 \varepsilon\Sigma_1 X_\varepsilon \Sigma_2^{-1} X_\varepsilon\Sigma_1X_\varepsilon}=I_d$.
This completes the proof.\hfill $\Box$

\bigskip
\begin{Remark}\label{AltForm}
The inverse of the solution of equation (\ref{Riccatti}) can be expressed in terms of $Y_\varepsilon$, the unique symmetric positive definite solution of the alternative
Riccati equation
\begin{equation}\label{AltRiccatti}
Y\Sigma_2 Y+\textstyle\frac \varepsilon 2 Y=\Sigma_1.
\end{equation}
In fact, if we write $Z_\varepsilon=\Sigma_2^{-1}X_\varepsilon \Sigma_1$ then
$$Z_\varepsilon=\Sigma_2^{-1}(\Sigma_2-\textstyle \frac \varepsilon 2 X_\varepsilon)X_\varepsilon^{-1}=(I_d-\textstyle \frac \varepsilon 2 \Sigma_2^{-1} X_\varepsilon)X_\varepsilon^{-1}=X_\varepsilon^{-1}-\frac \varepsilon 2 \Sigma_2^{-1}.$$
This shows that $Z_\varepsilon$ is symmetric. Also, since $\Sigma_2\geq \frac \varepsilon 2 X_\varepsilon$, we see that $X_\varepsilon^{-1}\geq\frac \varepsilon 2 \Sigma_2^{-1}$, that is, $Z_\varepsilon$ is positive definite. Since $Z_\varepsilon$ solves (\ref{AltRiccatti}) we conclude
$Z_\varepsilon=\Sigma_2^{-1}X_\varepsilon \Sigma_1=Y_\varepsilon$ or, equivalently, $X_\varepsilon \Sigma_1=\Sigma_2 Y_\varepsilon$.
From this we obtain $\Sigma_2^{-1}X_\varepsilon \Sigma_1X_\varepsilon=Y_\varepsilon X_\varepsilon$, which implies
\begin{eqnarray*}
I_d&=&\Sigma_2^{-1}(\Sigma_2-X_\varepsilon\Sigma_1 X_\varepsilon)+Y_\varepsilon X_\varepsilon=\textstyle\frac \varepsilon 2\Sigma_2^{-1}X_\varepsilon +Y_\varepsilon X_\varepsilon\\
&=&\big(\Sigma_2^{-1}+\textstyle\frac 2 \varepsilon Y_\varepsilon \big)\frac \varepsilon 2 X_\varepsilon.
\end{eqnarray*}
Thus, we conclude
\begin{equation}\label{AltFormEq}
\textstyle \frac 2 \varepsilon X_\varepsilon^{-1}=\Sigma_2^{-1}+\textstyle\frac 2 \varepsilon Y_\varepsilon.
\end{equation}

\end{Remark}

\bigskip
Before stating the announced result, we observe that in the analyisis of entropic regularization of transportation problems can 
focus on the case of centered probabilities $P$ and $Q$. In fact, for $\pi\in\Pi(P,Q)$ and $(X,Y)\sim \pi$ we write $\tilde{\pi}=\mathcal{L}(X-\mu_P,Y-\mu_Q)$,
$\tilde{P}=\mathcal{L}(X-\mu_P)$ and $\tilde{Q}=\mathcal{L}(Y-\mu_Q)$. The map $\pi\to \tilde{\pi}$ is a bijection between $\Pi(P,Q)$ and $\Pi(\tilde{P},\tilde{Q})$ and
$\int_{\mathbb{R}^d\times \mathbb{R}^d} \|x-y\|^2 d\pi(x,y)=\int_{\mathbb{R}^d\times \mathbb{R}^d} \|x-y\|^2 d\tilde{\pi}(x,y)+\|\mu_P-\mu_Q\|^2$. 
Similarly, we see that $H(\pi)=H(\tilde{\pi})$ and 
$K(\pi|\mu\otimes \nu)=K(\tilde{\pi}|\tilde{\mu}\otimes \tilde{\nu})$, where $d\tilde{\mu}(x)=d\mu(x-\mu_P)$, $d\tilde{\nu}(y)=d\nu(y-\mu_Q)$. If $\mu$, $\nu$ and
$\tilde{\mu}$, $\tilde{\nu}$ are equivalent,  we see, using (\ref{red1}), that
$$I_{\varepsilon,\mu,\nu}[\pi]=I_{\varepsilon,\mu,\nu}[\tilde{\pi}]+\|\mu_P-\mu_Q\|^2-\varepsilon\Big(\textstyle
\int_{\mathbb{R}^d} \log\big(\textstyle \frac {d\tilde{\mu}}{d\mu}(x)\big)d\tilde{P}(x)+\int_{\mathbb{R}^d} \log\big(\textstyle \frac {d\tilde{\nu}}{d\nu}(y)\big)d\tilde{Q}(y)
\Big).$$
With the choice of reference measures $\mu=\nu=N(0,\lambda I_d)$ we have $\tilde{\mu}=N(\mu_P,\lambda I_d)$,
$\tilde{\nu}=N(\mu_Q,\lambda I_d)$. Hence, $\log\big(\textstyle \frac {d\tilde{\mu}}{d\mu}(x)\big)=\frac{1}{2\lambda}(\|x\|^2-\|x-\mu_P\|^2)
=\frac{1}{2\lambda}(2\mu_P\cdot x-\|\mu_P\|^2)$ and we conclude that 
\begin{eqnarray}\label{W2ecent}
I_\varepsilon[\pi]&=&I_\varepsilon[\tilde{\pi}]+\|\mu_P-\mu_Q\|^2\\
\nonumber 
I_{\varepsilon,\lambda}[\pi]&=&I_{\varepsilon,\lambda}[\tilde{\pi}]+\|\mu_P-\mu_Q\|^2+\textstyle \frac {\varepsilon}{2\lambda} \Big(\|\mu_P\|^2+\|\mu_Q\|^2\Big).
\end{eqnarray}

\bigskip
\begin{Theorem}\label{GaussianEntropicTC}
If $P$ and $Q$ are Gaussian probabilities on $\mathbb{R}^d$ with means $\mu_1$ and $\mu_2$ and positive definite covariance
matrices  $\Sigma_1$ and $\Sigma_2$, respectively, then, if $\pi_0$ denotes the  Gaussian probability on $\mathbb{R}^d\times \mathbb{R}^d$
with mean $\mu=\big[ 
\begin{smallmatrix}
\mu_1 \\ \mu_2 
\end{smallmatrix}
\big]$ and covariance matrix $\Sigma_\varepsilon$ as in Proposition \ref{RiccattiSolutions}
\begin{eqnarray}
\nonumber
\mathcal{W}^2_{2,\varepsilon}(P,Q)&=& I_\varepsilon[\pi_0]=\|\mu_1-\mu_2\|^2\\
\label{entropiccost}
&&+\mbox{\em Tr}(\Sigma_1)+\mbox{\em Tr}(\Sigma_2)-2
\mbox{\em Tr}\big(\Sigma_1 X_\varepsilon \big){\textstyle  -\frac \varepsilon 2 \log\big((2\pi e)^{2d} (\frac{\varepsilon}2)^d|\Sigma_1 X_\varepsilon| \big) }.
\end{eqnarray}

\end{Theorem}

\noindent
\textbf{Proof.} {We write $r_P$ and $r_Q$ for the densities of $P$ and $Q$, respectively. From (\ref{W2ecent}) and the comments above we see that we only have to consider the case $\mu_1=\mu_2=0$. 
Also, since $H(\pi)$ can only be finite if $\pi$ has a density, we can rewrite (\ref{W2e}) as
$$\mathcal{W}^2_{2,\varepsilon}(P,Q)=\inf_{r\in \mathcal{R}(P,Q)}\Big[\int_{\mathbb{R}^d\times \mathbb{R}^d} [\|x-y\|^2 +\varepsilon \log r(x,y)] r(x,y)dx dy \Big]
$$
with $\mathcal{R}(P,Q)$ denoting the set of densities on $\mathbb{R}^d\times \mathbb{R}^d$ satisfying the marginal conditions
$\int_{\mathbb{R}^d} r(x,y)dy=r_P(x)$ for almost every $x$ and $\int_{\mathbb{R}^d} r(x,y)dx$ $=r_Q(y)$ for almost every $y$. Consider now $f\in L_1(P)$,
$g\in L_1(Q)$. Then for any $r\in \mathcal{R}(P,Q)$,
\begin{eqnarray*}
\lefteqn{\int [\|x-y\|^2 +\varepsilon \log r(x,y)] r(x,y)dxdy-\int f(x) dP(x)-\int g(y)dQ(y)}\hspace*{5cm}\\
&=&\varepsilon \int  r(x,y) \log\Big(\textstyle \frac{r(x,y)}{e^{\frac{f(x)+g(y)-\|x-y\|^2}\varepsilon}} \Big) dxdy\\
&\geq &
\varepsilon \int  e^{\frac{f(x)+g(y)-\|x-y\|^2}\varepsilon} \Big(\textstyle \frac{r(x,y)}{e^{\frac{f(x)+g(y)-\|x-y\|^2}\varepsilon}} -1\Big) dxdy\\
&=& \varepsilon -\varepsilon  \int  e^{\frac{f(x)+g(y)-\|x-y\|^2}\varepsilon}  dxdy,
\end{eqnarray*}
with equality if and only if $r(x,y)=e^{\frac{f(x)+g(y)-\|x-y\|^2}\varepsilon}$ for almost every $(x,y)$ (observe that this follows from the elementary fact
that $s\log s\geq s-1$, $s>0$, with equality if and only if $s=1$). This shows that
$$\mathcal{W}^2_{2,\varepsilon}(P,Q)\geq \varepsilon +\sup_{f\in L_1(P),g\in L_1(Q)} \big[\textstyle \int f(x) dP(x)+\int g(y)dQ(y)- 
\varepsilon  \int  e^{\frac{f(x)+g(y)-\|x-y\|^2}\varepsilon}  dxdy\big].$$
It shows also that if $r\in \mathcal{R}(P,Q)$ can be written as $r(x,y)=e^{\frac{f(x)+g(y)-\|x-y\|^2}\varepsilon}$ then $r$ is a minimizer for
the entropy-regularized transportation problem (indeed, by the strict convexity of $H$, the unique minimizer).

Now, if $\pi_0$ denotes the centered Gaussian distribution on $\mathbb{R}^d\times \mathbb{R}^d$ with 
covariance matrix $\Sigma_\varepsilon$ as in Proposition \ref{RiccattiSolutions}, then, obviously, $\pi_0\in \Pi(P,Q)$. From the expression
for $\Sigma_\varepsilon^{-1}$ and denoting $A_\varepsilon=\Sigma_1^{-1}+\frac 2 \varepsilon X_\varepsilon$ and $B_\varepsilon=\frac 2 \varepsilon X_\varepsilon^{-1}$ 
we see that the density of $\pi_0$ equals
\begin{eqnarray*}
r_0(x,y)&=&\textstyle \frac{1}{(2\pi)^d|\Sigma_\varepsilon|^{\frac 1 2}}\exp\Big[ 
-\frac 1 2 \big(x^TA_\varepsilon x+ y^T B_\varepsilon y-\textstyle \frac 4 \varepsilon x^T y\big)
\Big]\\
&=&
\textstyle \frac{1}{(2\pi)^d|\Sigma_\varepsilon|^{\frac 1 2}}\exp\Big[ 
-\frac 1 \varepsilon  \big(\|x-y\|^2+  x^T\big(\frac \varepsilon 2 A_\varepsilon -I_d  \big)x+ y^T\big(\frac \varepsilon 2 B_\varepsilon-I_d\big) y\big)
\Big].
\end{eqnarray*}
Consequently, $r_0(x,y)=e^{\frac{f_0(x)+g_0(y)-\|x-y\|^2}\varepsilon}$ with
\begin{eqnarray}\nonumber
f_0(x)&=& x^T\big(\textstyle I_d-\frac \varepsilon 2(\Sigma_1^{-1}+\frac 2 \varepsilon X_\varepsilon)\big)x-\frac \varepsilon 2 \log\big( (2\pi)^{2d} |\Sigma_\varepsilon|\big),\\
\label{g0}
g_0(y)&=& y^T\big(\textstyle I_d-X^{-1}_\varepsilon\big )y.
\end{eqnarray}
This proves that $\pi_0$ minimizes the regularized transportation cost between $P$ and $Q$.

Finally, to prove (\ref{entropiccost}) we note first that 
\begin{eqnarray}
\label{part1}
\int_{\mathbb{R}^d\times \mathbb{R}^d} \|x-y\|^2 d\pi_0(x,y)&=&
\mbox{Tr}(\Sigma_1)+\mbox{Tr}(\Sigma_2)-2
\mbox{Tr}(\Sigma_1 X_\varepsilon).
\end{eqnarray}
A simple computation shows that $H(\pi)=-\frac 1 2 \log\big((2\pi e)^{2d} |\Sigma_\varepsilon| \big)$. On the other hand
$$\det(\Sigma_\varepsilon)=\det(\Sigma_1) \det(\Sigma_2-X_\varepsilon \Sigma_1 \Sigma_1^{-1} \Sigma_1 X_\varepsilon)=
\textstyle\big(\frac{\varepsilon} 2\big)^d\det(\Sigma_1X_\varepsilon)$$ 
(here we have used that $\Sigma_2-X_\varepsilon \Sigma_1  X_\varepsilon=\frac \varepsilon 2 X_\varepsilon$).
Combining these last computations
with (\ref{part1}) we obtain (\ref{entropiccost}). \hfill $\Box$}

\bigskip
The proof of Theorem \ref{GaussianEntropicTC} can be easily adapted to other entropic regularizations. In particular, we can check that
$\pi_0$ is also the minimizer of $I_{\varepsilon,\lambda}[\pi]$ and also that
\begin{eqnarray}
\nonumber
\mathcal{W}^2_{2,\varepsilon,\lambda}(P,Q)&=& I_{\varepsilon,\lambda}[\pi_0]\\
\nonumber
&=&\|\mu_1-\mu_2\|^2+\textstyle \frac {\varepsilon}{2\lambda} \big(\|\mu_1\|^2+\|\mu_2\|^2\big)\\
&&\label{entropiccostlambda}
+\mbox{Tr}(\Sigma_1)+\mbox{Tr}(\Sigma_2)-2
\mbox{Tr}\big(\Sigma_1 X_\varepsilon \big)\\
&&\nonumber 
{\textstyle -\frac \varepsilon 2} \big[\log \big( |\Sigma_1 X_\varepsilon |\big) -{\textstyle \frac 1 \lambda \big(\mbox{Tr}(\Sigma_1)+\mbox{Tr}(\Sigma_2) \big)-d \big(2\log \lambda -\log {\textstyle \frac \varepsilon 2} -1\big)} \big].
\end{eqnarray}

\bigskip
Theorem \ref{GaussianEntropicTC} shows that the entropic transportation cost between normal laws is, as in the case of classical transportation cost,
a sum of two contributions. One accounts for the deviation in mean between the two laws. This part remains unchanged by the regularization with negative
differential entropy (but not with relative entropies). The other contribution, 
which accounts for deviations between the covariance matrices, behaves differently, but this behavior is also easier to understand in the case
of $\mathcal{W}_{2,\varepsilon}$. 
In the one-dimensional case we see that
$$\mathcal{W}^2_{2,\varepsilon}(N(0,\sigma_1^2),N(0,\sigma_2^2))=\sigma_1^2+\sigma_2^2-2\sqrt{\sigma_1^2\sigma_2^2+(\textstyle\frac{\varepsilon}4)^2}-\textstyle \frac \varepsilon 2 \log 
\big( \sqrt{\sigma_1^2\sigma_2^2+(\frac{\varepsilon}4)^2}-\frac \varepsilon 4\big)-\frac \varepsilon 2 \log ( 2 \pi^2 e \varepsilon ).$$ 
In particular, $\mathcal{W}^2_{2,\varepsilon}(N(0,1),N(0,1))=h(\frac \varepsilon 4)$ with
$$h(x)=2(1-\sqrt{1+x^2})-2x \log\big(\sqrt{1+x^2}-x\big)-2x\log \big((2\pi)^2 e x\big).$$ 
It is easy to see that $h(0)=0$, $h$ is decreasing in $\mathbb{R}_+$ and $\lim_{x\to\infty}h(x)=-\infty$.

\bigskip
While Theorem \ref{GaussianEntropicTC} is limited to Gaussian probabilities, its scope goes beyond that case. In classical optimal transportation
the Gaussian case provides a lower bound for the quadratic transportation cost through Gelbrich's bound (see in \cite{cuesta1996lower} which improves the bound in~\cite{gelbrich1990formula} ). We show next that this carries
over to entropic regularizations of transportation cost.

\begin{Theorem}\label{EntropicGelbrich}
If $P$ and $Q$ are probabilites on $\mathbb{R}^d$ with means $\mu_1,\mu_2$ and positive definite covariance matrices $\Sigma_1,\Sigma_2$, respectively,
then 
\begin{eqnarray}
\nonumber
\mathcal{W}^2_{2,\varepsilon}(P,Q)&\geq& \|\mu_1-\mu_2\|^2\\
\label{gelbound}
&&+\mbox{\em Tr}(\Sigma_1)+\mbox{\em Tr}(\Sigma_2)-2
\mbox{\em Tr}\big(\Sigma_1 X_\varepsilon \big){\textstyle  -\frac \varepsilon 2 \log\big((2\pi e)^{2d} (\frac{\varepsilon}2)^d|\Sigma_1 X_\varepsilon| \big) },
\end{eqnarray}
where $X_\varepsilon$ is as in (\ref{RiccattiSol}). Equality in (\ref{gelbound}) holds if and only if $P$ and $Q$ are Gaussian.
\end{Theorem}

\noindent
\textbf{Proof.} As in the proof of Theorem \ref{GaussianEntropicTC}, it suffices to consider the case of centered $P$ and $Q$. If $P$ (or $Q$)
does not have a density then $\Pi(P,Q)$ does not contain any probability with a density and, consequently, $I_\varepsilon[\pi]=+\infty$
for every $\pi\in \Pi(P,Q)$ and the result is trivial. We assume, therefore, that $P$ and $Q$ are absolutely continuous w.r.t. Lebesgue measure.
We consider $\pi\in \Pi(P,Q)$ with density $r$ and denote by $r_0$ the density of $\pi_0$, as defined in Theorem \ref{GaussianEntropicTC}. Then (recall (\ref{g0}))
\begin{eqnarray*}
I_\varepsilon[\pi]&=&\varepsilon \int_{\mathbb{R}^d\times\mathbb{R}^d} \log\Big({\textstyle \frac{r(x,y)}{e^{-\frac{\|x-y\|^2}{\varepsilon}}}}\Big) r(x,y) dx dy\\
&=&\varepsilon \int_{\mathbb{R}^d\times\mathbb{R}^d} \log\Big({\textstyle \frac{r(x,y)}{r_0(x,y)}}\Big) r(x,y) dx dy
+ \int_{\mathbb{R}^d\times\mathbb{R}^d} x^T (I_d-X_\varepsilon-{\textstyle\frac \varepsilon 2 \Sigma_1^{-1}})x r(x,y)dxdy\\
&&+ \int_{\mathbb{R}^d\times\mathbb{R}^d} y^T (I_d-X_\varepsilon^{-1})y r(x,y)dxdy-{\textstyle \frac \varepsilon 2 }\log \Big((2\pi)^{2d} ({\textstyle \frac \varepsilon 2})^d  |\Sigma_1 X_\varepsilon|\Big)\\
&=&\mbox{Tr}((I_d-X_\varepsilon-{\textstyle\frac \varepsilon 2 \Sigma_1^{-1}})\Sigma_1)+\mbox{Tr}((I_d-X_\varepsilon^{-1})\Sigma_2)-{\textstyle \frac \varepsilon 2 }\log \Big((2\pi)^{2d}({\textstyle \frac \varepsilon 2})^d|\Sigma_1 X_\varepsilon|\Big) +\varepsilon K(\pi|\pi_0),\\
&=&\mbox{Tr}(\Sigma_1)+\mbox{Tr}(\Sigma_2)-2\mbox{Tr}(\Sigma_1X_\varepsilon)-{\textstyle \frac \varepsilon 2 }\log \Big((2\pi e)^{2d}({\textstyle \frac \varepsilon 2})^d|\Sigma_1 X_\varepsilon|\Big) +\varepsilon K(\pi|\pi_0).
\end{eqnarray*}
Now (\ref{gelbound}) follows from the fact that $K(\pi|\pi_0)\geq 0$. If $P$ and $Q$ are Gaussian (and only in that case) then $\pi_0\in \Pi(P,Q)$. This completes the proof. \hfill $\Box$

\bigskip
To conclude this section we present a simple result on best approximation with respect to entropic transportation cost. In the case $\varepsilon=0$
(classical optimal transportation) $\mathcal{W}_2$ is a metric and for any $P$ with finite second moment we have
$$\mathcal{W}_2^2(P,Q)\geq  \mathcal{W}_2^2(P,P)=0, \quad Q\in \mathcal{F}_2(\mathbb{R}^d).$$
The fact that $\mathcal{W}_{2,\varepsilon}$ is no longer a metric for $\varepsilon>0$ changes the nature of the problem and we may wonder
which $Q\in \mathcal{F}_2(\mathbb{R}^d)$ is closest to $P$ in the sense of minimizing $\mathcal{W}_2^2(P,Q)$.
We show next that in the case of Gaussian $P$ the problem admits a simple solution.

\medskip
\begin{Theorem}\label{best_entropic_approximation}
Assume that $P$ is a probability on $\mathbb{R}^d$ with a density $r_P$ such that $\log r_P(x)\in L_1(P)$. Then   
$$P*N_d(0,{\textstyle \frac{\varepsilon} 2} I_d)=\mbox{\em argmin}_Q W_{2,\varepsilon}^2(P,Q),$$
with the minimization extended to the set of all probabilities on $\mathbb{R}^d$. Furthermore, $P*N_d(0,{\textstyle \frac{\varepsilon} 2} I_d)$
is the unique minimizer.
\end{Theorem}

\medskip
\noindent
\textbf{Proof.} We consider a probability on $\mathbb{R}^d\times \mathbb{R}^d$ with first marginal $P$ and 
$f\in L_1(P)$. Arguing as in the proof of Theorem \ref{GaussianEntropicTC} (take $g=0$) we see that 
$$I_\varepsilon[\pi]\geq  \varepsilon + \int f(x) dP(x) - \varepsilon \int  e^{\frac{f(x)-\|x-y\|^2}\varepsilon}  dxdy,$$
with equality if and only if $\pi$ has a density, $r$, that can be written as $r(x,y)=e^{\frac{f(x)-\|x-y\|^2}\varepsilon}$.
Now, if $r_0(x,y)=r_P(x) (\pi\varepsilon)^{-d/2} \exp(-\frac{\|y-x\|^2}{\varepsilon})$, then $r_0$ is a density on 
$\mathbb{R}^d\times\mathbb{R}^d$ with first marginal $P$, second marginal $P*N_d(0,\frac{\varepsilon} 2 I_d)$, and we can write
$r_0(x,y)=e^{\frac{f_0(x)-\|x-y\|^2}\varepsilon}$ with $f_0(x)=\varepsilon \log r_P(x)-\frac {d\varepsilon} 2 \log \pi\varepsilon$.
The assumption on $r_P$ ensures that $f_0\in L_1(P)$. We conclude that 
$$\min_Q W_{2,\varepsilon}^2(P,Q)= W_{2,\varepsilon}^2(P,P*N_d(0,{\textstyle \frac{\varepsilon} 2} I_d)).$$
Uniqueness follows by strict convexity of the entropic transportation cost.
\hfill $\Box$

\bigskip
We end this section with a simple observation that will be useful in our analysis of reguarized barycenters.
While $\mathcal{W}_2^2(P,Q)$ can take negative values, Theorem \ref{best_entropic_approximation} shows that the map $Q\mapsto \mathcal{W}_2^2(P,Q)$
is lower bounded by 
$$ W_{2,\varepsilon}^2(P,P*N_d(0,{\textstyle \frac{\varepsilon} 2} I_d))=\varepsilon \int_{\mathbb{R}^d}r_P(x) \log r_P(x)dx-{\textstyle \frac{d\varepsilon}{2}\log \pi\varepsilon}.$$

In the Gaussian case $P=N_d(\mu,\Sigma)$ we see that
$$ W_{2,\varepsilon}^2(N_d(\mu,\Sigma),Q)\geq W_{2,\varepsilon}^2(N_d(\mu,\Sigma),N_d(\mu,\Sigma+{\textstyle \frac \varepsilon 2} I_d))=
{\textstyle-\frac \varepsilon 2 \log |\Sigma| -\frac{d\varepsilon}2 \log (2\pi^2e\varepsilon)}.$$
This shows (recall Theorem \ref{GaussianEntropicTC}) that, in particular, if we fix a positive definite $\Sigma_1$ then the map 
$$\Sigma_2\mapsto  \mbox{Tr}(\Sigma_2)-2\mbox{Tr}(\Sigma_1X_\varepsilon) -\textstyle \frac{\varepsilon}{2}\log|\Sigma_1X_\varepsilon|,$$ 
with $X_\varepsilon$ as in (\ref{RiccattiSol}),
attains its minimal value within the set of positive definite matrices at $\Sigma_2=\Sigma_1+\frac \varepsilon 2 I_d$. 
Setting $A=\Sigma_1^{1/2}X_\varepsilon \Sigma_1^{1/2}$ is equivalent to setting $\Sigma_2=\Sigma_1^{-1/2}(A^2+\frac \varepsilon 2 A)\Sigma_1^{-1/2}$.
This allows to conclude that the strictly convex map (strict convexity follows easily from concavity of the log determinant)
$$A \mapsto  \mbox{Tr}(\Sigma_1^{-1}A^2)+ {\textstyle \frac \varepsilon 2}\mbox{Tr}(\Sigma_1^{-1}A)-2\mbox{Tr}(A) -\textstyle \frac{\varepsilon}{2}\log|A|$$ 
attains its minimal value within the set of positive definite matrices at $A=\Sigma_1$.

\section{Regularized barycenters.}

In this section we consider the entropic regularization of barycenters with respect to transportation cost metrics. To be precise, we will assume
that $P_1,\ldots,P_k$ are probabilities on $\mathbb{R}^d$ and $\lambda_1,\ldots,\lambda_k$ a collection of weights satysfying 
$\lambda_i>0$, $\sum_{i=1}^k \lambda_i=1$ and consider the functional
$$V_\varepsilon(Q)=\sum_{i=1}^k\lambda_i \mathcal{W}_{2,\varepsilon}^2(P_i,Q).$$
A minimizer of $V_\epsilon$ will be called an $\varepsilon$-regularized barycenter of $P_1,\ldots,P_k$ (with weights $\lambda_1,\ldots,\lambda_k$). 

The $L_2$ transportation cost metric, $\mathcal{W}_2$, satisfies the remarkable stability property that barycenters of Gaussian probabilities with
respect to $\mathcal{W}_2$ are Gaussian (this holds in fact for general location-scatter families, see  \cite{alvarez2016fixed}). We show in Theorem 
\ref{entropic_Gaussian_barycenter} that this
carries over to entropic regularized barycenters, although in this case the stability fails beyond the Gaussian case. 
Our result characterizes the barycenter in terms of the solution of a particular matrix equation, extending the result for the classical (unregularized) case. Existence and uniqueness of a solution for that matrix equation is guaranteed by our next result.

\medskip
\begin{Proposition} 
If $\Sigma_i\in\mathcal{M}_{d\times d}(\mathbb{R})$ are symmetric and positive definite then there exists a unique positive definite $\Sigma
\in\mathcal{M}_{d\times d}(\mathbb{R})$ such that
\begin{equation}\label{matrix_eq}
\sum_{i=1}^k \lambda_i \Big(\Sigma^{-1/2}\big( \Sigma^{1/2}\Sigma_i \Sigma^{1/2}+\textstyle (\frac{\varepsilon}{4})^2 I_d \big)^{1/2}\Sigma^{-1/2}+\textstyle \frac{\varepsilon}{4} \Sigma^{-1}\Big)=I_d.
\end{equation}
\end{Proposition}

\medskip
\noindent \textbf{Proof.}
The existence of a solution is equivalent to the existence of a fixed point for the map $G(\Sigma)=\sum_{i=1}^k\lambda_i G_i(\Sigma)$ with
$G_i(\Sigma)=\big( \Sigma^{1/2}\Sigma_i \Sigma^{1/2}+\textstyle (\frac{\varepsilon}{4})^2 I_d \big)^{1/2}+\textstyle \frac{\varepsilon}{4} I_d$.
$G$ is a continuous map on the set of positive semidefinite matrices and existence of a fixed point can be proved
using Brower's fixed point theorem, as follows. We write $A \preceq B$ to denote that $B-A$ is positive definite. Assume then that
$\alpha I_d\preceq \Sigma_i\preceq \beta I_d$, $i=1,\ldots,k$ and set $K=\{\Sigma :\, \alpha I_d\preceq \Sigma \preceq \beta+\frac \varepsilon 2 I_d\}$.
The set $K$ is compact and convex. Now, for every $\Sigma\in K$ we have $\Sigma^{1/2}\Sigma_i \Sigma^{1/2}+\textstyle (\frac{\varepsilon}{4})^2 I_d \preceq
\beta \Sigma + (\frac{\varepsilon}{4})^2 I_d\leq \big(\beta (\beta+\frac{\varepsilon} 2)+ (\frac{\varepsilon}{4})^2\big) I_d=
\big(\beta+\frac{\varepsilon}{4}\big)^2 I_d$. Using that $A\preceq B$ implies $A^{1/2}\preceq B^{1/2}$ (see, e.g., Theorem V.2.10 in \cite{Bhatia}) we conclude
that 
$$G_i(\Sigma)\preceq \textstyle{\big(\beta+\frac{\varepsilon}{2}\big)}I_d.$$
Similarly, we see that $\alpha I_d \preceq G_i(\Sigma)$ for $\Sigma\in K$. We conclude that $G$ maps $K$ into $K$ and from Brower's theorem we conclude
the existence of a fixed point. Uniqueness follows from Theorem \ref{entropic_Gaussian_barycenter} below \hfill $\Box$

\medskip
\begin{Theorem}\label{entropic_Gaussian_barycenter}
If $P_i=N(\mu_i,\Sigma_i)$ with $\mu_i\in \mathbb{R}^d$ and $\Sigma_i\in\mathcal{M}_{d\times d}(\mathbb{R})$ symmetric and positive definite then
the $\varepsilon$-regularized barycenter of $P_1,\ldots,P_k$ with weights $\lambda_1,\ldots,\lambda_k$ is $\bar{P}=N(\mu_0,\Sigma_0)$,
where $\mu_0=\sum_{i=1}^k \lambda_i \mu_i$ and $\Sigma_0$ is the unique positive definite solution of the matrix equation (\ref{matrix_eq}).
\end{Theorem}

%
%
%\medskip
%\noindent
%\textbf{Proof.} The comments before Theorem \ref{GaussianEntropicTC} and the fact that
%$\mu_0=\sum_{i=1}^k \eta_i \mu_i$ minimizes the function $m \mapsto \sum_{i=1}^k \eta_i \|\mu_i-m\|^2$ allow to consider only the case of centered $P_i=N(0,\Sigma_i)$.
%
%
%We note that for any
%$Q$ and any choice of $f_i\in L_1(P_i)$, $i=1,\ldots,k$ and $g_i\in L_i(Q)$, $i=1,\ldots,k$, satisfying $\sum_{i=1}^k\lambda_ig_i(y)=0$ for all $y$,
%$$V_\varepsilon(Q)\geq \varepsilon + \sum_{i=1}^k \lambda_i \Big[\int f_i(x_i) dP_i(x_i) - \varepsilon \int  e^{\frac{f_i(x_i)+g_i(y)-\|x_i-y\|^2}\varepsilon}  dx_idy\Big],$$
%with equality if and only if there are densities $r_1,\ldots,r_k$ with (common) second marginal equal to $Q$ such that
%$r_i(x_i,y)=e^{\frac{f_i(x_i)+g_i(y)-\|x_i-y\|^2}\varepsilon}$. A look at the proof of Theorems \ref{GaussianEntropicTC} and \ref{EntropicGelbrich}
%shows that equality can only happen for Gaussian $Q=N(0,\Sigma)$ such that 

%
%\bibliographystyle{plain}
%%\bibliographystyle{abbrv}
%\bibliography{thebib}

\end{document}